\documentclass[11pt]{amsart}
\usepackage{amsmath, amssymb, amsthm}
\usepackage{enumerate}
\newcommand{\taf}{{\hskip 5pt} $\blacksquare$
                  \renewcommand{\qedsymbol}{}}

\newcommand{\R}{{\mathbb{R}}}

\begin{document}
\title{The Fourier transform and the wave equation}
\author{Alberto Torchinsky}
\date{}
\maketitle

The study of  PDE's   arose in the 18th century in the context of
the development of models in the physics of continuous media, see \cite{BB}.
 It all began in 1747 when  d'Alembert, in a memoir
presented to the Berlin Academy,   introduced and analyzed the one
dimensional wave equation $u_{tt} =u_{xx} $ as a model of a
vibrating string, see \cite {dA},  and \cite{dA1}.  The
appropiatness of this model   still is  a topic of  discussion; for
recent work in this direction, see \cite{AK}, \cite {Bo}, and
\cite{LSC}.

  d'Alembert observed that a general solution to the wave equation is given by $u(x,t)=F(x+t)+G(x-t)$, where $F,G$
are arbitrary functions, and $F(x+t),G(x-t)$ represent   waves that
move  along the string to the left and right, respectively, with
constant speed.  The  Cauchy,  or initial value, problem for the one
dimensional wave equation consists of  finding   $u(x,t)$  in the
upper half-space $\R^2_+$ given its initial displacement $u(x,
0)=\varphi(x)$ and
 velocity $u_t(x, 0)=\psi(x) $. In terms of d'Alembert's traveling
 waves, the  solution
   assumes  the form
 \[u(x,t)= \frac{\varphi(x - t) + \varphi(x +  t)}2 + \frac12 \int_{x-t}^{x+t} \psi(y)\, dy\,,
\quad x\in\R, t>0\,.
\]
This expression serves as the prototype for the solution  to the
initial value problem for the  wave equation  in higher, odd,
dimensions.

The wave equation has a rich history  fraught with controversy,
both mathematical, see \cite{SL},
 physical, see \cite{WC},  and philosophical, see \cite{W}.
  At the heart of the   controversy   lies the notion of  a function,
as reflected by the  initial shape of the string. Contrary to
d'Alembert's view that only   $\varphi$'s and $\psi$'s that
could be described analytically should be allowed as    initial
data,
 Euler believed that  any curve, no
matter how irregular,  should  be included, see \cite {Eu}. At the same
time  D. Bernoulli, based on physical experiments,   conjectured
that the   shape of a vibrating string, i.e., $u(x,t)$, could be
described by means of a trigonometric series,  see  \cite {Be}, an
assertion that would be validated by Fourier's   work on the heat
equation, see \cite {F}.   Euler also believed  that partial
differential calculus should be extended to deal with irregular
functions, a notion  that was formalized in the  twentieth century
with the advent of generalized functions, or distributions, see \cite {Lu}, \cite {LS}.  These
two ideas: the use of Fourier transforms and their interpretation as
tempered distributions, are the essential ingredients in this paper.

Various   methods have been utilized since to solve the Cauchy
problem for the wave equation in higher dimensions - notably the
finite part of a divergent integral, see \cite {H}, the analytic
continuation of an integral of fractional order, see \cite {MR},
non-self-adjoint equations, see \cite{Cop}, and the Poisson
spherical means, see \cite {P}, \cite{FJ}, and \cite{FJ1} - see
\cite {Co} for a nice account of these results.

And, of course,   the Fourier transform. By means of the Fourier
transform in the space variables, a linear PDE with constant
coefficients is reduced to an ODE in the time variable, which, once
solved, gives the solution to the PDE in question  by means of the
Fourier inversion formula - the convolution properties of the
Fourier transform and the fact that the Fourier transform of a
spherically symmetric, or radial, function is again radial, are key
ingredients in this context. Of the three major examples of
classical second-order PDE's this approach provides a satisfactory
solution to the Laplace and heat equations, but,  for the wave equation
 the solution  lacks the clarity
of the other cases when the dimension $n$ is greater than $3$, see
\cite {CH}, \cite {Ta}, and \cite {T}.

In this paper I  will derive   the classical solution   to the wave
equation in dimensions $n$ greater than or equal to  $2$ via the
 Fourier analysis techniques outlined above,
 thus putting it on an equal footing with the other   equations.

Upon embarking on this endeavor it becomes quickly apparent that
  integrals of a function  on $\R^n$ that depends only on $x_n$  need to be
evaluated efficiently. Specifically,   the following formulas will
be used throughout the paper without further ado; for similar
results, see \cite{FJ1} and \cite{FS}.

{\it Suppose $n\ge 3$. If the function $f$ depends  on $x_n$ only,
\[\int_{ B(0,R)}f(x_n)\,dx=\omega_{n-1} \int_0^R\rho \,\int_{-\rho}^\rho f(s)\,\big(\rho^2-s^2\big)^{(n-3)/2}\,ds\,d\rho \,,
\]
and
\[\int_{\partial B(0,R)}f(x_n)\,d\sigma (x)=\omega_{n-1}\,R\int_{-R}^R f(s) \,(R^2-s^2)^{(n-3)/2}\,ds\,,\]
where $R>0$, $d\sigma$ is the element of surface area of the sphere
$\partial B(0,R)$, and $\omega_{n-1}$ is the surface area of the
unit ball of $\R^{n-1}$.}

To verify these formulas one  makes use of the $n$-dimensional polar
coordinates given by
\[\begin{cases}
x_1=\rho \, \sin(\varphi_1)\cdots \sin(\varphi_{n-3 })\, \sin(\varphi_{n-2 })\,\sin(\varphi_{n-1})\,,\\
x_2=\rho \, \sin(\varphi_1)\cdots \sin(\varphi_{n-3 })\,\sin(\varphi_{n-2 })\,\cos(\varphi_{n-1})\,,\\
x_3=\rho \, \sin(\varphi_1)\cdots \sin(\varphi_{n-3 })\,\cos(\varphi_{n-2 })\,,\\
\ldots\\
x_{n-1}=\rho\, \sin(\varphi_1)\,\cos(\varphi_2)\,,\\
x_n=\rho\, \cos(\varphi_1)\,,
\end{cases}
\]
where $0\le \varphi_k\le\pi$, $1\le k\le n-2$, $0\le
\varphi_{n-1}\le 2\pi$, and $0\le \rho$. As is readily seen by
induction on $n$, the Jacobian  of this transformation is
$J=\rho^{n-1}\sin^{n-2}(\varphi_1)\cdots
\sin^2(\varphi_{n-3})\,\sin(\varphi_{n-2})$.

Expressing the first volume integral above  in polar coordinates
gives
\[\int_{ B(0,R)}f(x_n)\,dx=c_n\,\int_0^R \rho^{n-1}\int_0^\pi f(\rho \cos(\varphi_1)\,)\,\sin^{n-2}(\varphi_1)\,d\varphi_1\,d\rho\,,\]
where $c_n$  is equal to
\[\int_0^\pi\ldots\int_0^\pi
\int_0^{2\pi} \sin^{n-3}(\varphi_2)\cdots \sin^2(\varphi_{n-3})\,\sin(\varphi_{n-2})\,
d\varphi_2\cdots d\varphi_{n-2}\,d\varphi_{n-1}\,.\]

Now, this expression is readily recognized as the surface area
$\omega_{n-1}$ of the unit sphere in $\R^{n-1}$, or, equivalently,
as the angular integrals corresponding to the evaluation   of the
volume $v_{n-1}$ of the unit ball in $\R^{n-1}$ in polar
coordinates, and so, $c_n=(n-1) \,v_{n-1} =\omega_{n-1}$.

As for the inner integral above,  note that it     equals
   \begin{align*}
\int_0^\pi f(\rho \cos(\varphi_1)\,)&\,\sin^{n-2}(\varphi_1)\,d\varphi_1\\
&=\int_0^\pi f(\rho \cos(\varphi_1)\,)\,(1-\cos^2(\varphi_1)\,)^{(n-3)/2}\,\sin(\varphi_1)\,d\varphi_1\,,
\end{align*}
which, with the change of variables $\rho \cos(\varphi_1)=s$,
becomes
\[\int_{-\rho}^\rho f(s)\,\big(1-(s/\rho)^2\big)^{(n-3)/2}\rho^{-1}\,ds
=\rho^{-(n-2)}\int_{-\rho}^\rho f(s)\,\big(\rho^2-s^2\big)^{(n-3)/2}\,ds\,.
\]

So, combining these observations it follows that
\[\int_{  B(0,R)}f(x_n)\,dx=
\omega_{n-1}\int_0^R\rho \,\int_{-\rho}^\rho
f(s)\,\big(\rho^2-s^2\big)^{(n-3)/2}\,ds\,d\rho \,,\] which is the
first relation  above. The second is obtained at once from this by
invoking the following   simple consequence of Fubini's theorem:
\[\int_{ B(0,R)}g(x)\,dx=\int_0^R\int_{\partial B(0,r)}g(x)\,d\sigma(x)\,dr\,.\]

Indeed, differentiating with respect to $R$ the identity
\[ \int_0^R\int_{\partial B(0,r)}f(x_n)\,d\sigma (x)\,dr=
\omega_{n-1}\int_0^R\rho \,\int_{-\rho}^\rho f(s)\,\big(\rho^2-s^2\big)^{(n-3)/2}\,ds\,d\rho\,, \]
it follows that
\[ \int_{\partial B(0,R)}f(x_n )\,d\sigma (x) =
\omega_{n-1}\, R \int_{-R}^R f(s)\,\big(R^2-s^2\big)^{(n-3)/2}\,ds\,,
\]
and the formula has been established.\taf

With this technical tool under the belt, the time has come to
 implement the Fourier analysis approach.  First, some notations.

  A multi-index $\alpha=(\alpha_1, \ldots,\alpha_n)$ is an $n$-tuple  of
non-negative integers.  The length $|\alpha|$ of $\alpha$ is
$|\alpha|=\alpha_1+\cdots+\alpha_n$, and, for
$x=(x_1,\dots,x_n)\in \R^n$, $x^\alpha= x_1^{\alpha_1}\cdots
x_n^{\alpha_n}$, and the differential monomial $D^\alpha$ is defined
as
\[D^\alpha \varphi (x)=\Big(\frac{\partial}{\partial x_1}\Big)^{\alpha_1}\ldots
\Big(\frac{\partial}{\partial x_n}\Big)^{\alpha_n}\,\varphi (x)\,.\]

Let $\mathcal S(\R^n)$ denote the Schwartz class of rapidly
decreasing functions, i.e., those
  functions $\varphi$ on $\R^n$ such that, together with their derivatives, decay faster than
any power of $x$ as $|x|\to\infty$. More precisely, $\varphi\in \mathcal S(\R^n)$ iff $\sup_{x\in\R^n}
\big| x^\alpha \,D^\beta\varphi(x)\big|<\infty$ for all multi-indices $\alpha,\beta$.

Recall that    the Fourier transform of a function $\varphi\in\mathcal S(\R^n)$ is given by
\[\widehat \varphi(\xi)=\int_{\R^n}\varphi(x)\,e^{-ix\cdot\xi}\,dx\,,
\]
  and that it also belongs to $\mathcal S(\R^n)$. In particular, note that
   $\widehat {D^\alpha \varphi}(\xi)= i^{|\alpha|}\xi^\alpha
  \widehat \varphi(\xi)$.
Furthermore,  the Fourier inversion formula holds and
\[  \varphi(x)
=\frac1{(2\pi)^n} \int_{\R^n}
\widehat \varphi(\xi)\,e^{i\xi\cdot x}\,d\xi\,.
\]

Also, if $\varphi*\psi(x)=\int_{\R^n} \varphi(x-y)\,\psi(y)\,dy$
denotes the convolution of $\varphi$ and $\psi$, it readily follows
that $\varphi*\psi\in\mathcal S(\R^n)$,
\[(\varphi*\psi)\, \widehat { } \  (\xi)=\widehat \varphi(\xi)\,\widehat \psi(\xi)\,,
\quad {\text {and, }}\ \varphi*\psi(x)=\frac1{(2\pi)^n} \int_{\R^n}
\widehat \varphi(\xi)\,\widehat \psi(\xi)\,e^{i\xi\cdot x}\,d\xi\,.
\]

 As for  the Fourier transform of a radial, or spherically symmetric, function $\varphi$ on $\R^n$, it is again radial.
Indeed, given a rotation matrix  $M$ in $\R^n$, since det$\,(M)=1$ and $ Mx\cdot M\xi=  x\cdot  \xi$, it follows that
 \[\widehat \varphi(M \xi)= \int_{\R^n}\varphi(x)\,e^{-ix\cdot M\xi}\,dx=
  \int_{\R^n}\varphi(M\mathbf{} x)\,e^{-i Mx\cdot M\xi}\,dx=
 \widehat \varphi(  \xi)\,.
 \]

 Now,
  if \[\nabla^2=\frac{\partial^2}{\partial x_1^2}+\cdots +  \frac{\partial^2}{\partial x_n^2}\]
denotes the Laplacian in $\R^n$,   the wave equation in the upper
half-space $\R^{n+1}_+$ is given by
\[\frac{\partial^2 u}{\partial t^2}=\nabla^2 u\,, \quad x\in \R^n, t>0\,.
\]

The Cauchy problem for this equation consists of finding
  $u(x,t)$ which satisfies the wave equation, subject to the initial conditions
  \[u(x,0)=\varphi(x)\,, \quad {\text{and}}\quad  u_t(x,0)=\psi(x)\,,\quad x\in\R^n ,
  \]
where,   partly because
of its simplicity,  we restrict ourselves to functions
 $\varphi,\psi$ in the Schwartz class $\mathcal S(\R^n)$.

Suppose $u$ solves this problem. Taking the Fourier transform in the
space variables, considering $t$ as a parameter, since
\[\widehat{u_{x_k x_k}}(\xi,t)=-\xi_k^2\, \widehat{u}(\xi,t)\,,\quad 1\le k\le n\,,\] it
readily follows that $\widehat {\nabla^2
u}(\xi,t)=-|\xi|^2\,\widehat u(\xi,t)$, and so, $\widehat{u}$
satisfies
\[ \widehat{u}_{tt}(\xi,t) +|\xi|^2 \widehat{u}(\xi,t)=0\,,\quad \xi\in \R^n, t>0\,,\]
subject to
\[ \widehat{u}(\xi,0)=\widehat{\varphi}(\xi),\quad {\text{and}}\quad
 \widehat{u_t}(\xi,0)=\widehat{\psi}(\xi)\,,\quad \xi\in\R^n \,.\]

 For each fixed $\xi\in\R^n$  this resulting ODE in $t$ is
 the simple harmonic oscillator equation with constant angular frequency $|\xi|$,
and one has
\[
\widehat{u}(\xi,t)=\widehat{\varphi}(\xi) \cos(t|\xi|) +
{\widehat{\psi}(\xi)} \,\frac{ \sin(t|\xi|)}{|\xi|}\,,\quad \xi\in \R^n, t>0\,.\]

At least formally, the inversion formula gives for
$(x,t)\in\R^{n+1}_+$,
\[
u(x,t)=\frac{1}{(2\pi)^n}\int_{\R^n} \widehat{\varphi}(\xi)
\cos(t|\xi|)\, e^{i\xi\cdot x}\,d\xi +\frac{1}{(2\pi)^n}\int_{\R^n}
{\widehat{\psi}(\xi)} \frac{ \sin(t|\xi|)}{|\xi|}\ e^{i\xi\cdot
x}\,d\xi\,.
\]

In fact, as long as $\widehat{\varphi} $ and $\widehat{\psi} $ are
sufficiently smooth so that   differentiation in $x$ and $t$ can be
performed under the integral sign, $u(x,t)$ will be a solution to
the Cauchy problem in the usual sense. However, a more transparent
expression that does not involve the Fourier transform of the
initial data is desirable, and, since the first integral above can
be formally obtained from the second by differentiating with respect
to $t$, I will concentrate on the latter.

Now, the  hope to interpret this integral as the inverse Fourier
transform of a  convolution
  involving   ${\psi} $ is dashed because  $\sin(|\xi|t)/|\xi| $
decays so slowly at infinity that it fails to be the Fourier transform of an integrable, or   square integrable,
function; this shortcoming is
  more apparent for $\cos(R|\xi|)$, which does not even tend
to $0$ as $|\xi|\to\infty$.

The next  best thing is then  to interpret   $\sin(|\xi|t)/|\xi| $ as
the Fourier transform of a  generalized function, or tempered
distribution. This will be a simple consequence of two identities
which, although classical in appearance,  are not familiar to many. The difference in the character of these identities when
the dimension is odd on the one hand and even on the other is not
unusual  in $n$-dimensional Fourier analysis; see \cite {SS} for
further illustrations of this phenomenon.

Specifically,   the identities in $\R^n$ are given by the following expressions.

{\it Assume that  $n$ is  an odd integer greater than or equal
to $3$. Then,
\[\frac{\sin (R|\xi|)}{|\xi|}= c_n\,\Big(\frac1{R}\,\frac{\partial}{\partial
R}\Big)^{(n-3)/2}
\Big( \frac{1}{\,\omega_n\,R}\int_{\partial B(0,R)}e^{-ix\cdot \xi}\,d\sigma(x)\Big)\,,
\]
where $R>0$, $\omega_n$ is the surface measure of the unit ball in
$\R^n$, and
$c_n^{-1}=(n-2)(n-4)\cdots 1 $.

On the other hand, if $n$ is an even integer greater than or equal
to $2$,
\[\frac{\sin(R|\xi|)}{|\xi|}= d_n\,\Big(\frac1{R}\,\frac{\partial}{\partial
R}\Big)^{(n-2)/2}
\Big( \frac{1}{\,v_n }\int_{  B(0,R)}\frac1{\big(R^2-|x|^2\big)^{1/2}}\
e^{-ix\cdot \xi}\,d x\Big)\,,
\]
where $R>0$,  $d_n^{-1}={n(n-2)(n-4)\cdots 2 }$,
and $v_n$ is the volume of the unit ball in $\R^n$.
}

Suppose first that $n$ is an odd integer greater than or equal to $
3$. One begins by computing $\widehat\chi_{ B(0,R)} (\xi)$, the
Fourier transform of the characteristic function of $B(0,R)$. Since
$\chi_{ B(0,R)}$ is radial, so is $\widehat\chi_{ B(0,R)}$, and,
consequently, $\widehat\chi_{ B(0,R)}(\xi)=\widehat\chi_{
B(0,R)}\big((0,\ldots, 0, -|\xi|)\big)$.
 Thus, the integral can be expressed as
\begin{align*}\int_0^R\int_{\partial  B(0,r)}e^{-ix\cdot
\xi}\,d\sigma (x)\,dr &= \widehat\chi_{ B(0,R)} (\xi)\\
&=\widehat\chi_{B(0,R)}\big((0,\ldots,0,-|\xi|)\big)\\
&=
 \int_{B(0, R)} e^{ix_n \,|\xi|}\,dx\\
&=\omega_{n-1}\int_0^R\rho \,\int_{-\rho}^\rho e^{is|\xi|}\,\big(\rho^2-s^2\big)^{(n-3)/2}\,ds\,d\rho \,.
\end{align*}
Hence, differentiating with respect to $R$,
\[\int_{\partial B(0,R)}e^{-ix\cdot \xi}\,d\sigma(x) =
\omega_{n-1} R \,\int_{-R}^R e^{is|\xi|}\,\big(R^2-s^2\big)^{(n-3)/2}\,ds\,,
\]
or, in other words,
\[
\int_{-R}^R e^{is|\xi|}\,\big(R^2-s^2\big)^{(n-3)/2}\,ds  =\frac1{\omega_{n-1} R }\,\int_{\partial B(0,R)}e^{-ix\cdot \xi}\,d\sigma (x)
\,.\]

Note that, since $n$ is an odd integer, $m=(n-3)/2$ is also an integer. Furthermore, since
$\big(R^2-s^2\big)^{m}$ vanishes at $R$ and $-R$, differentiating   with respect to $R$  gives
\[2mR\,\int_{-R}^R e^{is|\xi|}\,\big(R^2-s^2\big)^{m-1}\,ds=\frac{\partial}{\partial R}\Big(
\frac1{\omega_{n-1} R }\,\int_{\partial B(0,R)}e^{-ix\cdot \xi}\,d\sigma (x)\Big),
\]
or
\[\int_{-R}^R e^{is|\xi|}\,\big(R^2-s^2\big)^{m-1}\,ds=\frac1{2m}\,\Big(\frac1{R}\,\frac{\partial}{\partial R}\Big)
\Big( \frac1{\omega_{n-1} R }\,\int_{\partial B(0,R)}e^{-ix\cdot \xi}\,d\sigma (x)\Big).
\]

Clearly, repeating this process $m$ times yields
\[\int_{-R}^R e^{is|\xi|}\, ds=\frac1{2^m m!  }\,\Big(\frac1{R}\,\frac{\partial}{\partial R}\Big)^m
\Big(\frac1{\omega_{n-1} R }\,\int_{\partial B(0,R)}e^{-ix\cdot \xi}\,d\sigma (x)\Big).
\]

Therefore, since
\[\int_{-R}^R e^{is|\xi|}\,ds=\frac1{i|\xi|}\,\big(e^{iR|\xi|}-e^{-iR|\xi|}\big)=2\,\frac{\sin(R|\xi|)}{|\xi|}\,,
\]
it follows  that
\[\frac{\sin(R|\xi|)}{|\xi|}= c_n\,\Big(\frac1{R}\,\frac{\partial}{\partial R}\Big)^m
\Big( \frac{1}{\omega_n\,R}\int_{\partial B(0,R)}e^{-ix\cdot \xi}\,d\sigma(x)\Big)\,,
\]
 where $m=(n-3)/2$, and
 \[c_n=\frac{\omega_n}{\omega_{n-1}} \,\frac12\, \frac1{2^m m!}\,.\]

It only remains to compute $c_n$. Since
\[\frac{1}{\omega_n\,R}\int_{\partial
B(0,R)}d\sigma(x)=R^{n-2}\,,
\]
setting $|\xi|= 0$ in the relation just derived
 it follows that
\[R=
c_n\,\Big(\frac1{R}\,\frac{\partial}{\partial R}\Big)^m
R^{n-2}\,,
\]
and so, differentiating $m=(n-3)/2$ times obtains
\[R=(n-2)(n-4)\cdots 1\cdot c_n\,R\,,\]
and the conclusion follows in this case.

As for   even dimensions, the proof follows from the odd case by
essentially  Hadamard's    method of descent,  see \cite {H}.
  Since   $n$ is an even integer greater than or equal to $ 2$, $n+1$ is an odd integer
greater than or equal to  $ 3$. Observe that, for $\xi\in\R^n$,
$\xi'=(\xi,0)\in \R^{n+1}$, and
\[\frac{\sin(R|\xi|)}{|\xi|}= \frac{\sin(R|\xi'|)}{|\xi'|}\,,\quad \xi\in\R^n, R>0\,.\]

Now, for $x'\in\R^{n+1}$, since $x'=(x,x_{n+1})$ with $x\in\R^n$,
one has $x'\cdot \xi'=x\cdot \xi$. Moreover, in view that
$m'=(n+1-3)/2=(n-2)/2$, if $  B'(0,R)$ denotes the ball  in
$\R^{n+1}$ of radius $R$ centered at the origin, one has
\[\frac{\sin(R|\xi|)}{|\xi|}= c_{n+1}\,\Big(\frac1{R}\,\frac{\partial}{\partial R}\Big)^{m'}
\Big( \frac{1}{\omega_{n+1}\,R}\int_{\partial B'(0,R)}e^{-ix\cdot \xi}\,d\sigma(x')\Big)\,.
\]

Next  note that $\partial B'(0,R)$ can be decomposed into an  upper
and  a lower hemisphere. The upper hemisphere, $\partial B'(0,R)\cap
\{ x'_{n+1}\ge 0\}$,  is the graph of the function
$\gamma(x)=\big(R^2-|x |^2\big)^{1/2}$ for $x \in B (0,R)$, and,
similarly, the lower hemisphere, $\partial B'(0,R)\cap \{x'_{n+1}\le
0\}$, is the graph of $-\gamma$. Furthermore, since
\[\nabla \gamma(x)=-\frac1{\big(R^2-|x |^2\big)^{1/2}}\ x\,,\quad x\in\R^n,
\]
it follows that
\[d\sigma(x')= \sqrt{1+|\nabla \gamma(x)|^2}\,dx=
\frac{R}{\big(R^2-|x|^2\big)^{1/2}}\ dx\] in each hemisphere. Hence,
combining these observations, one has
 \[\frac {\sin(R|\xi|)}{|\xi|}=d_n
\,\Big(\frac1{R}\,\frac{\partial}{\partial R}\Big)^{(n-2)/2}
\Big( \frac1{v_n }\int_{  B (0,R)} \frac1{\big(R^2-|x|^2\big)^{1/2}}\,e^{-ix\cdot \xi}\,\ dx\Big)\,,
\]
where
\[ d_n =2\,\frac{c_{n+1}\, v_n}{\omega_{n+1}}
\,.\]

To evaluate $d_n$, we refer to the relation
\[c_{n+1}= \frac{\omega_{n+1}}{\omega_n}\, \frac1{2\,2^{m'}m'!}\]
established above. Since $n\,v_n=\omega_n$,
\[d_n=2\,\frac{c_{n+1}\, v_n}{\omega_{n+1}}= \frac1{n}\frac1{2^{m'}m'!}\,.
\]

Finally, in view that  $m'=(n-2)/2$,
\[2^{m'}m'! =2^{(n-2)/2}\big(\, (n-2)/2\big)!= 2^{(n-2)/2}\frac{1}{2^{(n-2)/2}}\, (n-2)(n-4)\cdots 2\,,\]
and, consequently,
$d_n \,{n(n-2)(n-4)\cdots 2}=1$.
\taf

Finally, all the pieces are in place to carry out the  program.
Suppose first that $n$ is an odd integer greater than or equal to
$3$. Two things must be proved: First, that the
   derivative  of
 the surface area measure carried on the sphere $\partial B(0,R)$ centered at the origin of
 radius $R$ in $\R^n$ given by
\[ T_{  R}=c_n\, \Big(\frac1{R}\,\frac{\partial}{\partial
R}\Big)^{(n-3)/2}
\Big( \frac{1}{\,\omega_n\,R}\int_{\partial B(0,R)} \,d\sigma(x)\Big)
\]
is a (compactly supported) tempered distribution, and, second, that
\[ \widehat{ T_{  R}}(\xi)= \frac{\sin(R|\xi|)}{|\xi|}\,,\quad \xi\in \R^n\,.
\]

Both of these assertions can be verified painlessly. Recall first
the notion of tempered distribution. A sequence
$\{\varphi_k\}\subset \mathcal S(\R^n)$ is said to tend to $0$ in
$\mathcal S(\R^n)$ provided that, for all multi-indices
$\alpha,\beta$, the sequence $\{x^\alpha D^\beta\varphi_k \}$
converges  uniformly to $0$  for every $x$ as $k\to\infty$, i.e.,
$\lim_{k\to\infty}\sup_x |x^\alpha D^\beta \varphi_k(x)|=0$. A
tempered distribution $T$ is a  linear functional $T:\mathcal
S(\R^n)\to \R$ that is continuous in the sense  that $T(
\varphi_k)\to 0$ whenever $\varphi_k\to 0$ in $\mathcal S(\R^n)$.

As for
 the Fourier transform  of a tempered distribution $T$,
 it  is defined as the tempered distribution $\widehat T$ that verifies
$\widehat T(\varphi)= T(\widehat \varphi)$, for every
$\varphi\in\mathcal S(\R^n)$;  that $\widehat T$ is a tempered
distribution is a consequence of the fact that,  if $\{\varphi_k\}$
tends to $0$ in $\mathcal S(\R^n)$,
 $\{\widehat{\varphi_k}\}$ also tends to $0$ in $\mathcal S(\R^n)$.
 Moreover, $T*\varphi$ is
defined for $\varphi\in\mathcal S(\R^n)$, say,   $(T*\varphi)\,
\widehat { } \ (\xi)=\widehat T(\xi)\,\widehat \varphi(\xi)$, and
the inversion formula holds in this case as well, see \cite{LS},
\cite{Ta}.

Now, to verify that $ T_{   R}$ is a tempered distribution observe
that, for $ \varphi\in\mathcal S(\R^n)$,
\[   T_{   R}(\varphi)   =c_n\, \Big(\frac1{R}\,\frac{\partial}{\partial
R}\Big)^{(n-3)/2}
\Big(R^{n-2}\,\frac{1}{\,\omega_n}\int_{\partial B(0,1)}\varphi(Rx')\,d\sigma(x')\Big).
\]

Since \[\Big(\frac1{R}\,\frac{\partial}{\partial
R}\Big)^{\ell}\varphi(Rx')\,,\quad 1\le \ell\le (n-3)/2\,,\] is a finite linear
combination of terms of the form $x'^\alpha D^\beta\varphi(Rx')$
with $|\alpha|,|\beta|\le \ell$ and coefficients that are
polynomials in $R$, a straightforward application of the Lebesgue
dominated convergence theorem (LDCT) gives that  the differentiation
may be performed under the integral sign. Furthermore, if
$\{\varphi_k\}$ tends to $0$ in $\mathcal S(\R^n)$,
$\lim_{k\to\infty}\sup_{x'}|x'^\alpha D^\beta\varphi(Rx')|=0$, and,
consequently,  again by the LDCT,
 $ T_{
  R}(\varphi_k)\to 0$ as $k\to\infty$. Hence, $ T_{  R}$ is
a (compactly supported) tempered distribution for each $R>0 $.

Next, there are a couple of ways of  computing $\widehat{T_{R}}$.
First, for $\varphi\in\mathcal S(\R^n)$,
\begin{align*}
T_R(\widehat\varphi)&= c_n \Big(\frac1{R}\,\frac{\partial}{\partial
R}\Big)^{(n-3)/2} \Big( \frac{1}{\,\omega_n\,R}\int_{\partial
B(0,R)}\widehat\varphi(x)\,d\sigma(x)\Big)
\\
&= c_n  \Big(\frac1{R}\,\frac{\partial}{\partial R}\Big)^{(n-3)/2}
\Big( \frac{1}{\,\omega_n\,R}\int_{\partial B(0,R)} \int_{\R^n}
  \varphi(\xi)\,e^{-i\xi\cdot x}\,d\xi\,d\sigma(x)\Big)\,.
\end{align*}

Now, since $  \varphi (\xi)e^{ -ix\cdot \xi}$ is integrable on
$\R^n\times  \partial B(0,R)$, by Fubini's theorem it is legitimate
 to change the order of integration,  and so,
\[ T_R(\widehat \varphi)=c_n\,
\Big(\frac1{R}\,\frac{\partial}{\partial R}\Big)^{(n-3)/2}\int_{\R^n}  \Big(
\frac{1}{\,\omega_n\,R} \int_{\partial B(0,R)}\,e^{- ix\cdot \xi} \,d\sigma(x)\Big)\,\varphi
(\xi)\,d\xi\,.
\]

A direct argument along the lines of the one used to
verify that $ T_{ R}$ is a tempered distribution gives that
 \[T_R(\widehat \varphi) =  \int_{\R^n}
 \frac{\sin(R|\xi|)}{|\xi|}\ \varphi (\xi)\,d\xi\,,
\]
and, consequently, $ \widehat {T_R}= {\sin(R|\xi|)}/{|\xi|}$ in the
sense of distributions.

On the other hand, one may note  that, since $T_{R}$ has compact
support, by elementary properties of the Fourier transform of
tempered distributions, $\widehat{T_{R }}$ is
 the $C^\infty(\R^n)$ function given by $\widehat{
T_{R}}(\xi)=T_R( e^{-i x \cdot \xi})$, and, consequently, again by the
 identity involving $\sin(R|\xi|)/|\xi|$, $\widehat{ T_{R}}(\xi) = \sin(R|\xi|)/|\xi|$.

Therefore, the integral  of interest to us can be written as
\[
\frac{1}{(2\pi)^n}\int_{\R^n} {\widehat{\psi}(\xi)} \frac{
\sin(t|\xi|)}{|\xi|}\ e^{i\xi\cdot x}\,d\xi
=\frac{1}{(2\pi)^n}\int_{\R^n} {\widehat{\psi}(\xi)}
\widehat{T_t}(\xi)\,e^{i\xi\cdot x}\,d\xi\,,\] and so, by the
inversion formula  for tempered distributions, this  expression is
equal to $T_t*\psi(x)$. So, it only remains to identify $T_t*\psi(x)$.

Since by definition $T*\psi$ is the $ C^\infty(\R^n)$ function
$T*\psi(x)=\tau_x T(\widetilde \psi)$, where $\tau_x$ denotes
translation by $x$ and $\widetilde \psi(y)=\psi(-y)$, it readily
follows that
\[T_t * \psi(x)=c_n\,
\Big(\frac1{t}\,\frac{\partial}{\partial
t}\Big)^{(n-3)/2}
\Big( \frac{1}{\,\omega_n\,t} \int_{\partial B(x,t)}\, \psi(y)\,d\sigma(y)\Big)\,.
\]

In other words, if
  \[M_R  \varphi (x)=\frac{1}{\omega_n R^{n-1}}\int_{\partial B(x,R)} \psi(y)\,d\sigma(y)
  \] denotes  the average   of
  $\psi$ over the sphere centered at  $x$ of radius $R$,
 \[T_t * \psi(x)=c_n\,
\Big(\frac1{t}\,\frac{\partial}{\partial
t}\Big)^{(n-3)/2}
\big(  t^{n-2} \, M_t   \psi (x) \big)\,.
 \]

So, the solution to the Cauchy problem in this case can be expressed compactly
in terms of spherical means as follows, see \cite {E}.

{\it
Let $n$ be an odd integer greater than or equal to $3$, and suppose
that $\varphi,\psi$ are smooth functions. Then, with
$c_n^{-1}=(n-2)(n-4)\cdots 1 $,
 \begin{align*} u(x,t)=  {c_n}\,\frac{\partial}{\partial t}\,\Big(\frac1{t}\frac{\partial}{\partial
t}  \Big)^{(n-3)/2}
 &\Big( t^{n-2}M_t \varphi (x)\Big)\\
&+
 {c_n}\, \Big(\frac1{t}\frac{\partial}{\partial
t}\Big)^{(n-3)/2}
\Big( t^{n-2}M_t \psi (x)\Big)
\end{align*}
is a $C^2(\R^{n+1}_+)$ function that satisfies $u_{tt}(x,t)=\nabla^2
u(x,t)$, $ x\in \R^n,\, t>0$, and $u(x,0)=\varphi(x)$, $
u_t(x,0)=\psi(x)$, $x\in \R^n$.

In fact, when $u$ is defined as above in terms of spherical means,
it suffices to assume that $\varphi\in C^{(n+3)/2}(\R^n)$
and $\psi\in C^{(n+1)/2}(\R^n)$.
 }

When $n=3$   the expression reduces      to the familiar Kirchhoff's formula.

The  result  for even dimensions follows mutatis mutandis from the
odd dimensional case. Indeed, if $n$ is an even integer greater than
or equal to $ 2$ and $d_n^{-1}=n\,(n-2)\cdots 2$, let
\[T_t  =d_n
\,\Big(\frac1{t}\,\frac{\partial}{\partial t}\Big)^{(n-2)/2}
\Big( \frac1{v_n }\int_{  B (0,t)} \frac1{\big(t^2-|x|^2\big)^{1/2}} \,\ dx\Big)\,.
\]

 Then, $T$ is a compactly supported tempered distribution,  and $\widehat{T_t}(\xi)= \sin(t|\xi|)/|\xi|$. Thus, if
  \[\mathcal M_R  \varphi(x)=\frac{1}{v_n R^{n}}\int_{B(x,R)}\frac1{\big(R^2-|x-y|^2\big)^{1/2}} \ \psi(y)\,dy
  \] denotes  the weighted average   of
  $\psi$ over the ball centered at  $x$ of radius $R$,
 \[T_t * \psi(x)=d_n\,
\Big(\frac1{t}\,\frac{\partial}{\partial
t}\Big)^{(n-2)/2}
\big(  t^n  \, \mathcal M_t   \psi (x) \big)\,.
 \]

Therefore, the following holds, see \cite {E}.

{\it
Let $n$ be an even integer greater than or equal to $2$, and suppose
that $\varphi,\psi$ are smooth functions. Then, with
$d_n^{-1}=n\,(n-2)\cdots 2$,
 \begin{align*} u(x,t)= {d_n}\,\frac{\partial}{\partial t}\,\Big(\frac1{t}\frac{\partial}{\partial
t}  \Big)^{(n-2)/2}
 &\Big( t^{n }\mathcal M_t \varphi (x)\Big)\\
&+
d_n\, \Big(\frac1{t}\frac{\partial}{\partial
t}\Big)^{(n-2)/2}
\Big( t^{n }\mathcal M_t \psi (x)\Big)
\end{align*}
is a $C^2(\R^{n+1}_+)$ function that satisfies $u_{tt}(x,t)=\nabla^2
u(x,t)$, $ x\in \R^n,\, t>0$, and $u(x,0)=\varphi(x)$, $
u_t(x,0)=\psi(x)$, $x\in \R^n$.

In fact, when $u$ is defined as above in terms of   means,
it suffices to assume that $\varphi\in C^{(n+4)/2}(\R^n)$
and $\psi\in C^{(n+2)/2}(\R^n)$.
 }

In   particular, when $n=2$, the result   is known as  Poisson's formula.

I am grateful to the referee and the editor for their valuable
comments that shaped the final presentation of this paper.

\end{document}